\input amstex
\documentstyle{amsppt}
\magnification=1200
\UseAMSsymbols
\NoBlackBoxes
\topmatter
\title structure of total subspaces of dual banach spaces
\endtitle
\rightheadtext{TOTAL SUBSPACES OF DUAL BANACH SPACES}
\author M.I.Ostrovskii
\endauthor
\address Mathematical Division, Institute for Low Temperature Physics
and Engineering, 47 Lenin avenue, 310164 Kharkov, UKRAINE \endaddress
\email ostrovskii\%ilt.kharkov.ua\@relay.ussr.eu.net \endemail
\thanks It is an English translation of the paper published
in Teoriya Funktsii, Funktsional'nyi Analiz i ikh Prilozheniya, vol. 58
(1992), 60--69\endthanks
\endtopmatter
I. Let $X$ be a Banach space, $X^*$ its dual. The unit ball and  the
unit sphere of $X$ are denoted by $B(X)$ and $S(X)$ respectively.  Let  us
recall some definitions.

A subspace $M$ of $X^*$ is said to be {\it total} if for every 
$0\neq x\in X$ there
is an $f\in M$ such that $f(x)\neq 0$. A subspace $M\subset X^*$ 
is said to be {\it norming} if
for some $c>0$ we have
$$
(\forall x\in X)(\sup _{f\in S(M)}|f(x)|\ge c||x||).
$$
Let $M$ be a subspace of $X^*$. The  set  of  all  limits  of  weak$^*$
convergent sequences from $M$ is called its {\it weak}$^*$  
{\it sequential  closure}
and is denoted by $M_{(1)}$. It is clear that $M_{(1)}$ is a  linear  subspace
of $X^*$. But $M_{(1)}$ need not be closed and, consequently,  need  not  be
weak$^*$ closed. Corresponding example has been given by S.Mazurkiewicz.
This was a reason to S.Banach to introduce [1, pp.~208, 213]  weak$^*$
sequential closures (S.Banach used  the  term  "derive  faible")  of
other orders, including transfinite ones. For ordinal $\alpha $  
the  {\it weak}$^*$
{\it sequential closure of order} $\alpha $ of a subspace 
$M\subset X^*$ is the set 
$$M_{(\alpha )}=
\bigcup _{\beta <\alpha }(M_{(\beta )})_{(1)}.
$$
For a chain of weak$^*$ sequential closures we have
$$
M_{(1)}\subset M_{(2)}\subset \ldots 
\subset M_{(\alpha )}\subset M_{(\alpha +1)}\subset \ldots .
$$
If we have $M_{(\alpha )}=M_{(\alpha +1)}$ then all  subsequent  closures  
coincide
with $M_{(\alpha )}$. The least ordinal $\alpha $ for which 
$M_{(\alpha )}=M_{(\alpha +1)}$ is  called  the
{\it order} of the subspace M.

The study of total subspaces with infinite orders turns out  to
be important in the theory of  topological  vector  spaces \cite{4},\cite{9},
subspaces of order 2  turns  out  to  be  useful  in  the  theory  of
improperly posed problems \cite{11}. So the  problem  of  description  of
total subspaces of prescribed order arises in  a  natural  way.  The
purpose of the present paper is to investigate the following version
of this problem. Let a separable Banach space $X$ and an ordinal $\alpha $ are
given. For what Banach spaces $Y$ there exists an isomorphic embedding
$T:Y\to X^*$ such that $T(Y)$ is a total subspace of order $\alpha $ in $X^*$?

Let us say few words about terminology and notation.  The  term
"operator" means a bounded linear operator. For  a  subset $A$ of  a
Banach space $X$ lin$(A)$ and cl$(A)$ are, respectively, the  linear  span
of $A$ and the closure of $A$ in the strong topology. For a subset $A$ of
a dual Banach space $X^*$ $w^*-$cl$(A)$  and $A^\top$  are,  respectively,  the
closure of $A$ in the weak$^*$ topology and the set 
$\{x\in X:(\forall x^*\in A)(x^*(x)=0)\}$.
The direct sum of spaces $X$ and $Y$ is denoted by $X\oplus Y$. 
We hope that our
terminology and notation  are  standard  and  self-explanatory.  Our
sources for Banach space basic concepts and results are \cite{7},\cite{8},
\cite{14}.

Now we shall list known results about weak$^*$ sequential closures
which will be used in this paper. Let $X$ be a separable Banach space.

1. Subspace $M\subset X^*$ satisfies the equality $M_{(1)}=X^*$ if and only if
$M$ is norming \cite{1, p.~213}.

2. If the Banach space $X$ is quasireflexive  then  \cite{14, p.~78}
every total subspace $M$ of $X^*$ is norming and hence we have 
$M_{(1)}=X^*$.

3. For any subspace $M\subset X^*$ we have \cite{1, p.~124}:
$$
(M=M_{(1)})\Rightarrow (M=w^*-cl(M)).
$$

4. The order of any subspace of $X^*$ is a countable \cite{14, p.~50}
non-limit \cite{6} ordinal.

5. If the Banach space $X$ is a nonquasireflexive  then  for  any
countable ordinal $\alpha $ the space $X^*$ contains a total subspace of order
$\alpha +1$, \cite{10}.

II. Let us consider the problem stated in the paragraph I.  The
case when $X$ is reflexive is trivial. If $X$ is quasireflexive then  by
statement 2 from I the order of any total subspace of $X^*$ is 0 or 1.
It is easy to see that in this case the set of isomorphic  types  of
total subspaces of order 1 coincides  with  the  set  of  isomorphic
types of closed subspaces of $X^*$  with  codimension  between  1  and
$\dim (X^{**}/X)$. The question whether this subspaces are isomorphic to $X^*$
is the version of well-known problem. We  shall  not  consider  this
question.

If the space $Y$ is such that $Y^*$ does not contain closed  norming
subspaces of infinite codimension and $Y$ is  isomorphic  to  a  total
subspace of $X^*$ then by \cite{12, Theorem 3.1} $Y$  is  isomorphic  to  a
norming subspace of $X^*$ and hence $X$ is isomorphic to  a  subspace  of
finite codimension of $Y^*$.
On the other hand it is easy to see that  if  the  space $Y$  is
nonquasireflexive and $X$  is  isomorphic  to  a  finite-codimensional
subspace of $Y^*$ then $Y$ is isomorphic to a total  (and  even  norming)
subspace of $X^*$.
Therefore if $Y$ is nonquasireflexive space such that $Y^*$  doesn't
contain closed norming subspaces of infinite  codimension  then  the
operator $T$ from the problem posed in  the  first  paragraph  may  be
found only for $\alpha =0$ or $\alpha =1$ and if and only if $X$ 
is isomorphic  to $Y^*$
or to finite-codimensional subspace of $Y^*$ respectively.

For quasireflexive space all that we say  above  remains  true.
The only additional condition is that $X$  must  be  isomorphic  to  a
subspace of $Y^*$ with codimension no greater than $\dim (Y^{**}/Y)$.
Therefore  it  is  natural  to  consider  only  the   following
particular cases of the problem from the first paragraph.

Let $X$ be a nonquasireflexive separable Banach space. Let $Y$ be a
Banach space isomorphic to  a  subspace  of $X^*$  and  such  that $Y^*$
contains closed norming subspaces of infinite codimension. Let $\alpha $  be
a countable ordinal.

{\bf Question 1.} Does there exist an isomorphic embedding $T:Y\to X^*$ for
which the subspace $T(Y)$ is total?

{\bf Question 2.} For what countable ordinals  does  there  exist  an
isomorphic embedding $T:Y\to X^*$ such that 
$(T(Y))_{(\alpha )}\neq (T(Y))_{(\alpha +1)}=X^*$?

At first we shall show  that  in  general  the  answer  to  the
question 1 is negative. After this (in the third paragraph) we shall
find additional condition under which the answer onto the question 1
is positive and consider question 2.
\proclaim{Theorem 1} There  exist  a  separable  Banach  space $X$  and  a
separable  subspace $Y\subset X^*$  such  that $Y^*$  contains  closed  norming
subspaces of infinite codimension but $Y$ is not isomorphic to a total
subspace of $X^*$.
\endproclaim

We need space constructed in \cite{2}. Let $1<p<2.$ Let us  denote  by
$X_n (n=0,1,2,\ldots)$ a subspace of $l_p$ consisting of  the  vectors  whose
coordinates are equal to zero begining with $(n+1)$-th.
The space $J(l_p)$ is the completion of the space of all  finitely
non-zero sequences $\{x_n\}^{\infty }_{n=0} (x_n\in X_n)$ under the norm
$$
||\{x_n\}^{\infty }_{n=0}||_J=\sup (\sum^{m-1}_{k=1}
||x_{p(k+1)}-x_{p(k)}||^2+||x_{p(m)}||^2)^{1/2},
$$
where the $\sup $ is over all  increasing  integer  
sequences $\{p(i)\}^m_{i=1}$
with $p(1)\ge 0.$

We shall use the following result about the spaces $J(l_p), 1<p<2.$
\proclaim{Theorem 2} \cite{13} 
\roster
\item"A." The space $(J(l_p))^{**}$ can be represented as a
direct sum $Z\oplus l_p$, where $Z$ is the canonical  image  
of $J(l_p)$  in  its
second dual.
\item"B." Every weakly null bounded away  from  0  sequence  in $J(l_p)$
contains a subsequence equivalent to the orthonormal basis of $l_2$.
\endroster
\endproclaim
Let us pass to the proof of theorem 1. Let $X=l_q\oplus (J(l_p))^*$, where
$1/q=1-1/p$. Let $Y$ be a canonical image of $J(l_p)$ in $X^*$.  Theorem 2.A
implies  immediately  that $Y^*$  contains  closed  separable  norming
subspaces of infinite codimension (e.g. the annihilator in $Y^*$ of the
direct summand isomorphic to $l_p$ in $Y^{**})$.
It remains to check that $Y$  is  not  isomorphic  to  any  total
subspace of $X^*$. In order to do this we need the following  variation
on the theme of Pitt's result \cite{8, p.~76}.
\proclaim{Lemma 1} Every operator $R:J(l_p)\to l_p (1<p<2)$ is compact.
\endproclaim
\demo{Proof} Assume that cl$(R(B(J(l_p))))$ is a non-compact  set.  Then
there is a sequence $\{y_i\}^{\infty }_{i=1}$ in $B(J(l_p))$ so that
$$
(\exists \delta >0)
(\forall i,j\in {\Bbb N}, i\neq j)(||Ry_i-Ry_j||\ge \delta ).
$$
Since $(J(l_p))^*$ is separable then we can find  a  weakly  Cauchy
subsequence $\{y_{i(n)}\}^{\infty }_{n=1}$ in 
$\{y_i\}^{\infty }_{i=1}$. For a sequence $z_n=y_{i(2n)}-y_{i(2n-1)}\
(n\in {\Bbb N})$ we have
$$
(\forall n\in {\Bbb N})(||Rz_n||\ge \delta );
\eqno{(1)}$$
$$
w-\lim_{n\to \infty }z_n=0.
\eqno{(2)}$$
By theorem 2.B we can select in $\{z_n\}^{\infty }_{n=1}$ a subsequence 
$\{z_{n(k)}\}^{\infty }_{k=1}$
equivalent to the orthonormal basis of $l_2$. By (1), (2) and well-known
results  \cite{8, p.~7, 53}  the  sequence $\{Rz_{n(k)}\}^{\infty }_{k=1}$   
contains   a
subsequence equivalent to the unit vector basis of $l_p$. Since $p<2$  we
arrived at a contradiction. The lemma is proved.
\enddemo

Let $T:Y\to X^*$  be  arbitrary  isomorphic   embedding.   We   have
$X^*=l_p\oplus (J(l_p))^{**}=l_p\oplus J(l_p)\oplus l_p$. 
Let us denote by $P_1, P_2$  and $P_3$  the
projections   corresponding   to   this   decomposition.   We   have
$T=P_1T+P_2T+P_3T$. The operators $P_1T$ and $P_3T$ are  compact  by Lemma  1.
Hence the kernel of $P_2T$ is finite-dimensional. Let us represent  the
space $Y$ as a direct sum $\ker P_2T\oplus Y_1$. The compactness of 
operators $P_1T$
and $P_3T$ implies that the restriction of $P_2T$ to $Y_1$ is an isomorphism.
We shall denote it by $R$. The subspace 
$R(Y_1)\subset \{0\}\oplus J(l_p)\oplus \{0\}$  will  be
denoted by $M$. Subspace $T(Y_1)\subset X^*$ can be 
represented in the  following
form:
$$
T(Y_1)=\{(P_1TR^{-1}m, m, 
P_3TR^{-1}m):m\in M\}\subset l_p\oplus J(l_p)\oplus l_p.
$$
It is clear that in order to finish the proof of Theorem  1  it
is sufficient to prove that weak$^*$ closure of $T(Y)$  in $X^*$  does  not
coincide with $X^*$. Since the linear span of the union of weak$^*$ closed
and  finite-dimensional  subspaces  is  weak$^*$  closed  then  it   is
sufficient to prove that $w^*-$cl$(T(Y_1))$ is of infinite codimension  in
$X^*$.

Since $M\subset \{0\}\oplus J(l_p)\oplus \{0\}$ then every operator 
$Q:M\to X^*$ is  continuous
if  we endow both spaces with the topology $\sigma (X^*,X)$. Let us introduce
operator $Q:M\to X^*$ defined in the following way:
$$
Q(0,m,0)=(P_1TR^{-1}m,m,P_3TR^{-1}m).
$$
It follows from what we have said above that $Q$ is an isomorphic
embedding. It is easy to see that $Q$ has (unique) $\sigma (X^*,X)$-continuous
extension  onto  the  subspace $w^*-$cl$(M)\subset X^*$.  
Let  us  denote   this
extension by $Q_e$. It is easy to check that $Q_e$ can be represented as a
sum of the identity operator (on $w^*-$cl$(M)$) and a  compact  operator.
Therefore the kernel of $Q_e$ is finite-dimensional and the restriction
of $Q_e$ to any weak$^*$ closed complement of ker$Q_e$ in $w^*-$cl$(M)$ is a
$\sigma (X^*,X)$-continuous isomorphism. Therefore by  the  Krein--Smulian
theorem the image of $Q_e$ is weak$^*$ closed. Since $Q(M)=T(Y_1)$ then we
obtain $w^*-$cl$(T(Y_1))\subset $im$Q_e$. 
On the other hand since $P_1Q$ is  a  compact
operator and $w^*-$cl$(M)$ is contained in 
$\{0\}\oplus J(l_p)\oplus l_p$ then $P_1Q_e$ is also
compact. Therefore im$Q_e$ and hence also $w^*-$cl$(T(Y_1))$ is of infinite
codimension in $X^*$. The theorem is proved.

III. In this part of the paper additional condition under which
answer onto Question 1 become affirmative is pointed out.  We  shall
also consider Question 2.

Let $X$ be a Banach space. Let $Y$  be  a  subspace  of $X^*$.  Every
element of $X$ may  be  considered  in  a  natural  way  as  a  linear
functional on $Y$. So there is a natural mapping of $X$ into $Y^*$. We shall
denote this mapping by $H_Y$.

\proclaim{Theorem 3}
\roster
\item"A." Let $X$ be a separable Banach space.  Let $Y$  be  a
subspace of $Y^*$ such that $H_Y(X)$ is of  infinite  codimension  in $Y^*$.
Then there exists an isomorphic embedding $T:Y\to X^*$ the image of which
is total and, moreover, we have $(T(Y))_{(2)}=X^*$.
\item"B." If additionally the subspace $Y^\top\subset X$  is  
infinite-dimensional
then  the  operator $T$  can  be  chosen   in   such   a   way   that
$(T(Y))_{(1)}\neq (T(Y))_{(2)}=X^*$.
\item"C." If additionally the subspace $Y^\top\subset X$ is nonquasireflexive  
then
for every countable ordinal $\alpha $ there exists an isomorphic embedding
$T:Y\to X^*$ such that $(T(Y))_{(\alpha )}\neq (T(Y))_{(\alpha +1)}=X^*$.
\endroster
\endproclaim
\remark{Remark} Condition of the part A of Theorem 3 are  satisfied  if
the image of $H_Y$ is non-closed.
\endremark
\demo{Proof of Theorem  3} Using  the  arguments  from  the  proof  of
Proposition 3.5 of [12] we can prove the following result.
\proclaim{Lemma 2} Let $X$ be a Banach space. Let $Y$ be  a  subspace  of $X^*$
such that $H_Y(X)$ is  non-closed.  Then  there  exists  an  isomorphic
embedding $E:Y\to X^*$ such that cl$(E^*(X))$ is of infinite  codimension  in
$Y^*$ and $(E(Y))^\top\supset Y^\top$.
\endproclaim

By this lemma we  may  restrict  ourselves  to  the  case  when
cl$(H_Y(X))$ is of infinite codimension in $Y^*$.
By separability of $X$ and remark from \cite{3, p.~358} there exists a
weak$^*$ null basic sequence $\{u_i\}^{\infty }_{i=1}$ in $Y$ such that 
for certain bounded
sequence $\{v^{**}_k\}^{\infty }_{k=1}\subset X^{**}$ and 
certain partition $\{I_k\}^{\infty }_{k=1}$ of the positive
integers into pairwise disjoint infinite subsets we have
$$v^{**}_k(u^*_i)=\cases 1,&\text{if $i\in I_k$;}\\
0,&\text{if $i\not\in I_k$.}\endcases$$

At first we shall prove parts A and B of Theorem  3.  Let
$\{x^*_k\}^{\infty }_{k=1}$ be a normalized sequence spanning a 
norming subspace of $X^*$.
Let us introduce the operator $T:X^*\to X^*$ given by
$$
T(x^*)=x^*+\sum^{\infty }_{k=1}4^{-k}v^{**}_k(x^*)x^*_k/||v^{**}_k||.
\eqno{(3)}$$
It is easy to see that for every $x^*\in X^*$ we have
$$
(2/3)||x^*||\le ||T(x^*)||\le (4/3)||x^*||.
\eqno{(4)}$$
Therefore $T$ is an isomorphism. Let us show  that $(T(Y))_{(2)}=X^*$.
For $i\in I_k$ we have $T(u^*_i)=u^*_i+4^{-k}x^*_k/||v^{**}_k||$. 
Since the sequence $\{u^*_i\}^{\infty }_{i=1}$ is
weak$^*$ null we obtain that $x^*_k\in (T(Y))_{(1)}$. 
Since vectors $\{x^*_k\}^{\infty }_{k=1}$ span
a norming subspace in $X^*$ it follows that $(T(Y))_{(2)}=X^*$.  Part  A  is
proved.

In order to prove part B we have to prove that in the case when
$Y^\top\subset X$ is in\-fi\-ni\-te-di\-men\-sio\-nal the subspace 
$T(Y)\subset X^*$ is nonnorming. Let
$y\in Y$ and $z\in (\cap ^n_{k=1}\ker  x^*_k)\cap Y^\top$. 
Then 
$$(Ty)(z)=\sum^{\infty }_{k=n+1}4^{-k}v^{**}_k(y)x^*_k(z)/||v^{**}_k||.$$
Therefore $|(Ty)(z)|\le (4^{-n}/3)||y||||z||$.     
By     (4)     we     have
$$|(Ty)(z)|\le (4^{-n}/2)||Ty||||z||.$$ 
Hence subspace $T(Y)\subset X^*$ is nonnorming. Part
B is proved.

Let us turn to the proof of part  C.  Since  the  space $Y^\top$  is
nonquasireflexive then by result  5  from  the  first  part  of  the
present  paper  there  exists  a  subspace $M$  of $(Y^\top)^*$  such  that
$$M_{(\alpha )}\neq M_{(\alpha +1)}=(Y^\top)^*$$ 
if ordinal $\alpha $ is infinite and such that 
$$M_{(\alpha -1)}\neq M_{(\alpha )}=(Y^\top)^*$$ 
if ordinal $\alpha \ge 1$ is finite. Let $N$ be the set of all bounded
extensions of functionals from $M$ onto the whole $X$. It is proved  in
\cite{10} that 
$$N_{(\alpha )}\neq N_{(\alpha +1)}=X^*\ (N_{(\alpha -1)}\neq 
N_{(\alpha )}=X^*).$$  
By  separability  of $X$
there exists a normalized sequence 
$\{x^*_k\}^{\infty }_{k=1}$ such that $x^*_k\in N (k\in {\Bbb N})$  and
for $L$=lin$(\{x^*_k\}^{\infty }_{k=1})$ we have 
$L_{(\alpha )}\neq L_{(\alpha +1)}=X^*$  if $\alpha $  is  infinite  and
$L_{(\alpha -1)}\neq L_{(\alpha )}=X^*$ if $\alpha $ is finite.

Let  us  introduce  operator $T:X^*\to X^*$  by  equality  (3).  This
operator satisfies inequality (4) and hence is an  isomorphism.  Let
us show that the subspace $(T(Y))\subset X^*$ satisfy the condition
$$
(T(Y))_{(\alpha )}\neq (T(Y))_{(\alpha +1)}=X^*.
$$
At first let us show that
$$
(T(Y))_{(1)}\subset \hbox{cl(lin}(Y_{(1)}\cup \{x^*_k\}^{\infty }_{k=1})).
$$
Let $z\in (T(Y))_{(1)}$ and let the sequence 
$\{z_i\}^{\infty }_{i=1}\subset T(Y)$ be such  that
$z=w^*-\lim_{i\to \infty }z_i$. We have $z_i=T(y_i)$ 
for certain vectors $y_i\in $Y. Since $T$ is
an isomorphism then the sequence $\{y_i\}^{\infty }_{i=1}$ is bounded. 
We may consider
without loss of generality that this sequence is  weak$^*$  convergent.
Operator $T$ is a sum of identical and compact operators.  Denote  the
second one by $K$. We may assume without loss of generality  that  the
sequence $\{K(y_i)\}^{\infty }_{i=1}$ is strongly convergent. 
It  is  clear  that  its
limit belongs to cl(lin$(\{x^*_k\}^{\infty }_{k=1}))$.  Hence, 
$z\in$ cl(lin$(Y_{(1)}\cup \{x^*_k\}^{\infty }_{k=1}))$.
Therefore restrictions of functionals from $(T(Y))_{(1)}$ to the subspace
$(Y^\top)\subset X$ belong to $M$. This implies that 
$(T(Y))_{(\alpha )}\subset N_{(\alpha -1)}$ in the case
when ordinal $\alpha $ is  finite,  and 
$(T(Y))_{(\alpha )}\subset N_{(\alpha )}$  in  the  case  when
ordinal $\alpha $ is infinite. In both cases we have 
$(T(Y))_{(\alpha )}\neq X^*$.

Using the same arguments as in the proof of part A we can show
that $x^*_k\in (T(Y))_{(1)}$. Therefore $L\subset (T(Y))_{(1)}$. 
Hence $(T(Y))_{(\alpha +1)}=X^*$ both
when $\alpha $ is finite and when $\alpha $ is infinite. 
The proof of Theorem 3 is complete.
\enddemo
\remark{Remarks}
\roster
\item"1." Let $\alpha \ge 2$ be a non-limit ordinal. By arguments  which
were used in the proofs of Theorem 3 and Theorem  3.1  of \cite{12}  it
follows that Banach space $Y$ is isomorphic to a subspace of  order $\alpha $
in the dual of some separable Banach space if and only if the  space
$Y^*$  contains  closed  separable   norming   subspace   of   infinite
codimension.
\item"2." If spaces $X$ and $Y$ satisfy the conditions of the  part  A  of
Theorem 3 and $X$  is  isomorphic  to  a  direct  sum $X\oplus Z$  for  some
in\-fi\-ni\-te-di\-men\-sio\-nal (nonquasireflexive) 
Banach space $Z$, then  there
exists an isomorphic embedding $Q:Y\to X^*$ the image of  which  satisfies
the conditions of the part B (C) of Theorem 3.
\item"3." If $X$ is a separable Banach space, $Y$ is a subspace of $X^*$  and
$Y^*$ is nonseparable then $X$ and $Y$ satisfy the conditions of part A  of
Theorem 3. By \cite{7, p.~147, 213} it is the case if $X$  is  a  separable
${\Cal L}_{\infty }$-space (in the sense of  Lindenstrauss  -  Pelczynski)  
and $Y$  is nonreflexive subspace of $X^*$.
\item"4." Let $X=C(K)$ where $K$ is a metrizable compact. It follows  from
well-known properties of the spaces  of  continuous  functions  (see
\cite{7}) and Remarks 2 and 3 that for any nonreflexive subspace 
$Y\subset X^*$ and
every countable ordinal $\alpha \ge 1$  there  exists  an  isomorphic  
embedding
$T:Y\to X^*$  such  that $(T(Y))_{(\alpha )}\neq (T(Y))_{(\alpha +1)}=X^*$.  
It  is  easy  to  see
immediately that the same is true for $\alpha =0.$
\endroster
\endremark
\proclaim{Theorem 4} Let $X=l_1\oplus Z$, where $Z$ is  a  separable  
Banach  space.
Then for any countable ordinal $\alpha \ge 1$ and any Banach space $Y$  
the  dual
of which contains closed separable  norming  subspaces  of  infinite
codimension there exists an isomorphic embedding $T:Y\to X^*$ such that
$(T(Y))_{(\alpha )}\neq (T(Y))_{(\alpha +1)}=X^*$. 
Conversely, the dual of any total subspace
$Y\subset X^*$  contains  closed  separable  norming  subspaces  of   infinite
codimension.
\endproclaim
\demo{Proof} Let $a,b>0.$ Recall that a  subset $A\subset X^*$  is  
said  to  be
$(a,b)$-{\it norming} if the following conditions are satisfied:
$$
(\forall x\in X)(\sup \{|x^*(x)|:x^*\in A\}\ge a||x||);
$$
$$
\sup \{||x^*||:x^*\in A\}\le b.
$$
Let Banach spaces $X$ and $Y$ be such that there exists an operator
$R:X\to Y^*$ for which the set $R(B(X))$ is $(a,b)$-norming subset of 
$Y^*$  for
some $a>0$ and $b<+\infty $, and subspace $R(X)$ is of infinite  
codimension  in
$Y^*$. It is easy to see that in this case the restriction of $R^*$ to the
subspace $Y\subset Y^{**}$ is  an  isomorphic  embedding,  and $X$  and 
$R^*(Y)\subset X^*$
satisfy the conditions of part A of Theorem 3.
By the quotient-universality of $l_1$\cite{8, p.~108} all what is said
above is valid for $X=l_1\oplus Z$ and any Banach space $Y$, 
the dual of  which
contains closed separable norming subspaces of infinite codimension.
Using Remark 2 and part C of Theorem 3 we obtain the  proof  of  the
first part of Theorem 4.

Let us turn to the proof of the second part. If  subspace $Y\subset X^*$
is nonseparable then subspace cl$(H_Y(X))\subset Y^*$  is  a  closed  separable
norming subspace of infinite codimension.

Turn to the case when $Y$ is separable. If subspace $Y$ is not
norming then repeating arguments from the proof of  Proposition  3.5
of \cite{12} we find in $Y^*$ closed separable norming subspace of  infinite
codimension.

If subspace $Y$ is norming then we finish the proof with the use
of the following generalization of Lemma 1 of \cite{5}.

\proclaim{Lemma 3} Let $Y$ be a separable Banach space, whose dual contains
a subspace isomorphic to $l_1$. Then $Y^*$ contains  norming  subspace  of
infinite codimension isomorphic to $l_1$.
\endproclaim
\demo{Proof} Let $\{f^*_i\}^{\infty }_{i=1}$ be a sequence in $Y^*$ which is  
equivalent  to
the unit vector basis of $l_1$. The space $Y$ is separable and  so  there
exists a weak$^*$ convergent subsequence $\{f^*_{i(n)}\}^{\infty }_{n=1}$. 
Then the sequence
$g^*_n=f^*_{i(2n)}-f^*_{i(2n-1)}\ (n\in {\Bbb N})$  is  a  weak$^*$  
null  sequence   which   is
equivalent to the unit vector basis of $l_1$. Let $0<c,C<\infty $ be such that
$$
(\forall \{a_i\}^{\infty }_{i=1}\subset {\Bbb R})
(c\sum |a_i|\le ||\sum a_ig^*_i||\le C\sum |a_i|).
$$
Let $\{y_i\}^{\infty }_{i=1}$ be dense in $S(Y)$. 
For every $i\in {\Bbb N}$ we choose a functional
$h^*_i\in S(Y^*)$ such that $h^*_i(y_i)=1.$ Let $\{m(i)\}^{\infty }_{i=1}$ 
be an increasing sequence
of even positive integers satistying the condition:
$$
(\forall i\in {\Bbb N})(|g^*_{m(i)}(y_i)|<c/4).
$$
It  is  clear  that  the  sequence $z^*_i=(c/2)h^*_i+g^*_{m(i)}\ 
(i\in {\Bbb N})$  is
equivalent to the unit vector basis of $l_1$. Let 
$L=$cl(lin$\{z^*_i\}^{\infty }_{i=1})$.  It
remains  to  prove  that $L$  is  a  norming  subspace  of   infinite
codimension of $Y^*$.
We have $||z^*_i||\le C+(c/2)$ and $|z^*_i(y_i)|>(c/2)|h^*_i(y_i)|-
|g^*_{m(i)}(y_i)|>c/4$.
Since $\{y_i: i\in {\Bbb N}\}$ is a dense subset  of 
$S(Y)$  then $L$  is  a  norming
subspace of $Y^*$.
We shall show that
$$
\hbox{lin}(\{g^*_{2k-1}\}^{\infty }_{k=1})\cap L=\{0\},
$$
from which it follows that $L$ is of infinite codimension in $Y^*$.  Let
us suppose that sequences $\{b_k\}^{\infty }_{k=1}$ and 
$\{a_i\}^{\infty }_{i=1}$ from $l_1$  satisfy  the
equality
$$
\sum^{\infty }_{k=1}b_kg^*_{2k-1}+\sum^{\infty }_{i=1}a_ig^*_{m(i)}=
(c/2)\sum^{\infty }_{i=1}a_ih^*_i.
$$
The norm of the right-hand side of this equality is not greater
than $(c/2)\sum |a_i|$, and the norm of the left-hand side is not less than
$c(\sum |b_k|+\sum |a_i|)$. Hence all $a_i$  and  all 
$b_k$  are  equal  to  0.  This
completes the proof of the lemma.
\enddemo
\enddemo

\enddocument